\documentclass[a4paper,11pt]{amsart}

\usepackage{microtype}
\usepackage{times}
\usepackage{geometry}
\usepackage{amsmath}
\usepackage{amsfonts}
\usepackage{amsthm}
\usepackage{amssymb}
\usepackage{upref}
\usepackage{color}
\usepackage{graphicx}
\usepackage{hyperref}
\usepackage{array}

\graphicspath{{./anisop_art/}{../figures/}{../figures_eps/}}

\theoremstyle{plain}
\newtheorem{thm}{Theorem}[section]

\newtheorem{prop}[thm]{Proposition}
\newtheorem{deff}[thm]{Definition}
\newtheorem{rem}[thm]{Remark}
\theoremstyle{definition}
\newcommand{\ds}{\displaystyle}
\newcommand{\Per}{\operatorname{Per}}
\newcommand{\di}{\operatorname{div}}
\newcommand{\mrestr}{%
  \,\raisebox{-.127ex}{\reflectbox{\rotatebox[origin=br]{-90}{$\lnot$}}}\,%
}
\newcommand{\bo}[1]{{\bf #1}}
\newcommand{\gconv}{\stackrel{\Gamma}{\longrightarrow}}
\newcommand{\wconv}{\rightharpoonup}

\title{Numerical study of optimal partitioning problems in relation to an anisotropic perimeter}
\author{Beniamin Bogosel}

\begin{document}

\begin{abstract}
We present a $\Gamma$-convergence approximation for the total anisotropic length of a partition. This
theoretical result gives rise to a numerical method which allows the study of minimal partitions with respect
to different anisotropies. We also give a numerical framework for the study of isoperimetric problems with
density.
\end{abstract}

\maketitle

\section{Introduction} 

If we consider a set $\Omega \subset \Bbb{R}^d$ with $C^1$ boundary, then its perimeter is equal to
\[ \Per(\Omega)=\int_{\partial \Omega} 1 d\mathcal{H}^{n-1} =\int_{\partial \Omega} \|\vec n(x)\| d\mathcal{H}^{n-1} \]
where $\vec n(x)$ denotes the unit outer normal vector corresponding to $x \in \partial \Omega$. Thus, the perimeter treats all directions in the same way and no direction has an advantage over the others. Things change if we pick another norm $\varphi$ on $\Bbb{R}^d$, different from the euclidean one. We can define the anisotropic perimeter associated to a norm $\varphi$ by
\[ \Per_\varphi(\Omega)=\int_{\partial \Omega} \varphi(\vec n) .\]

The problem of studying the partitions which minimize total perimeter of the cells has been studied before. A famous result due to Hales \cite{hales} states that every partition of the plane into sets of equal areas has perimeter greater than the hexagonal honeycomb tiling. In $\Bbb{R}^3$ the problem of finding the optimal tilling with respect to the perimeter using shapes of equal volume is still open. Kelvin conjectured that truncated octahedra may be optimal, but Weaire and Phelan \cite{weaire-phelan} found a better configuration. 

The study of the partitions which minimize the sum of anisotropic perimeters is even more challenging, since the optimal partition depends on the norm $\varphi$. This motivates the interest in providing efficient numerical algorithms which compute the optimal partitions. One such method was developed by \'E. Oudet in \cite{oudet} in the isotropic case in two and three dimensions. The author uses $\Gamma$-convergence to approximate the sum of perimeters of the parts by a relaxation using $\Gamma$-convergence and a theorem of Modica and Mortola. A different approach to minimal partitions in an anisotropic setting with applications to image classification is presented in \cite{amstutz-aniso}.

The first main contribution of this article is to prove a general $\Gamma$-convergence result regarding the approximation of the sum of anisotropic perimeters of a partition. The anisotropic variant of Modica-Mortola's theorem can be found, for example in \cite{braides2}, \cite{braides} or \cite{ambrosio-aniso}. The framework we present also allows that the anisotropy $\varphi$ depends also on the position of the point. As a particular case, we obtain also an approximation by $\Gamma$-convergence for the density perimeter given by $\int_{\partial \Omega} \nu(x)d \mathcal{H}^{n-1}(x)$. 

Secondly, we provide a numerical method which is efficient in the study of partitions minimizing the total anisotropic perimeter. We give a few examples where we favorize a fixed number of directions. 
We are also able to study numerically isoperimetric problems concerning a density perimeter. In particular we obtain some results proved in \cite{camivi} and \cite{morgan-pratelli}.

\section{Preliminaries}

The notion of $\Gamma$-convergence, introduced by de Giorgi, is a suitable tool for the study of the convergence of variational problems. Its properties make it a suitable tool for the study of limits of variational problems.

\begin{deff}
Let $X$ be a metric space and $F_\varepsilon, F: X \to [0,+\infty]$ a sequence of functionals on $X$ (defined for $\varepsilon>0$). We say that $F_\varepsilon$ $\Gamma$-converges to $F$ and we denote $F_\varepsilon \gconv F$ if the following two properties hold:

\begin{itemize}
\item[(LI)] For every $x \in X$ and every $(x_\varepsilon) \subset X$ with $x_\varepsilon \to x$ we have
\begin{equation} F(x) \leq \liminf_{\varepsilon \to 0} F_\varepsilon(x_\varepsilon)\label{LI} \end{equation}

\item[(LS)] For every $x \in X$ there exists $(x_\varepsilon) \subset X$ such that $(x_\varepsilon) \to x$ and
\begin{equation} F(x) \geq \limsup_{\varepsilon \to 0} F_\varepsilon(x_\varepsilon).\label{LS}\end{equation}
\end{itemize}
\label{def-g-conv}
\end{deff}
Given $x_0 \in X$ we will call \emph{recovery sequence} a sequence $(x_\varepsilon)$, which satisfies property \eqref{LS}. This sequence satisfies, in particular, the relation
\[ \lim_{\varepsilon \to 0} F_\varepsilon(x_\varepsilon)=F(x).\]
Here are three main properties of the $\Gamma$-convergence. 

\begin{prop}
If $F_\varepsilon \gconv F$ in $X$ then the following properties hold:
\begin{itemize}
\item[(i)] $F$ is lower semicontinuous;
\item[(ii)] If $G: X \to [0,\infty)$ is a continuous functional then 
\[ F_\varepsilon+G \gconv F+G.\]
\item[(iii)] Suppose $x_\varepsilon$ minimizes $F_\varepsilon$ over $X$. Then every limit point of $(x_\varepsilon)$ is a minimizer for $F$.
\end{itemize}
\label{properties-gconv}
\end{prop}


Sometimes it is difficult to prove the (LS) property \eqref{LS} for every $x \in X$. Having an element $x$ with some good regularity properties may aid in constructing the recovery sequence. The following procedure, of reducing the class of elements $x$ for which we prove \eqref{LS} to a dense subset of $\{F<+\infty\}$, is classical (see for example \cite{braides2},\cite{braides}).

\begin{prop}
 Let $\mathcal{D} \subset  \{F<+\infty\}$ be a dense subset of $X$, such that for every $x \in \{F<+\infty\}$ and $(u_n) \subset \mathcal{D}$, with $(u_n) \to x$ we have
\[ \limsup_{n \to \infty} F(u_n) \leq F(x). \]
 Suppose that for every $x \in \mathcal{D}$, the  property \eqref{LS} is verified. Then \eqref{LS} is verified in general. 
\label{reduction-regular}
\end{prop}

\begin{rem} In general if $F_\varepsilon \stackrel{\Gamma}{\longrightarrow} F$ and $G_\varepsilon \stackrel{\Gamma}{\longrightarrow} G$ we cannot conclude that $ F_\varepsilon + G_\varepsilon \stackrel{\Gamma}{\longrightarrow} F+G.$ 
Thus, the result proved in Section \ref{main-results} is not trivial. One sufficient condition for the above implication to hold would be that for each $u$ we could find the same recovery sequence for $F$ and $G$. For more details and examples see \cite{braides2}.
\end{rem}

The notion of $\Gamma$-convergence was introduced in Definition \ref{def-g-conv} and its main properties were stated in Proposition \ref{properties-gconv}. One classical $\Gamma$-convergence result is the Modica Mortola theorem. For the sake of completeness, we rewrite its statement below. For simplicity, we denote 
\[ X = \{ u \in L^1(D) : \int_D u = a\},\]
where $a \in (0,|D|)$ is a fixed constant. 

\begin{thm} (Modica-Mortola) Let $D$ be a bounded open set and let $W: \Bbb{R} \to [0,\infty)$ be a continuous function such that $W(z)=0$ if and only if $z \in \{0,1\}$. Denote $c =2\int_0^1 \sqrt{W(s)}ds$. We define $F_\varepsilon,F : L^1(D) \to [0,+\infty]$ by
\[ F_\varepsilon(u) = \begin{cases}\varepsilon \int_D|\nabla u|^2 +\frac{1}{\varepsilon} \int_D W(u) & u \in H^1(D) \cap X\\ +\infty & \text{otherwise} \end{cases} \]
and
\[ F(u)             = \begin{cases} c\Per(u^{-1}(1)) & u \in BV(D; \{0,1\})\cap X \\ +\infty & \text{otherwise}\end{cases}\]
then
\[ F_\varepsilon \stackrel{\Gamma}{\longrightarrow} F\]
in the $L^1(D)$ topology.
\label{modica-mortola2}
\end{thm}
For a proof of this result we refer to \cite{gammaconvalberti} or \cite{buttazzogconv}. In the numerical simulations we usually choose the potential \[W(s)=s^2(1-s)^2\] which gives the corresponding constant $c=1/3$. The numerical importance of this theorem was just recently observed. Indeed, when one wants to compute numerically the perimeter of a set $\Omega$, the boundary $\partial \Omega$ must be well known. Using a parametric formulation might work if one only has to deal with one set. As soon as we consider multiple shapes which might touch, keeping track of each parametrized boundary is not a simple task. If we want to study a partitioning problem, using a parametric formulation rises difficulties in imposing the non-overlapping condition. This is a point where having a good relaxation for the perimeter, like the theorem mentioned above, becomes really useful.

\begin{figure}
\includegraphics[width=0.19\textwidth]{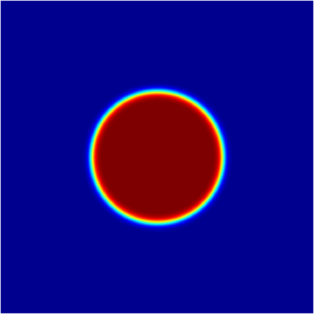}~
\includegraphics[width=0.19\textwidth]{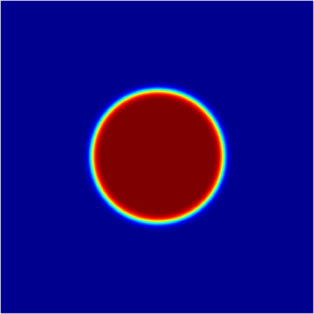}~
\includegraphics[width=0.19\textwidth]{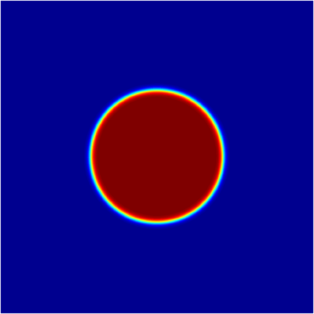}~
\includegraphics[width=0.19\textwidth]{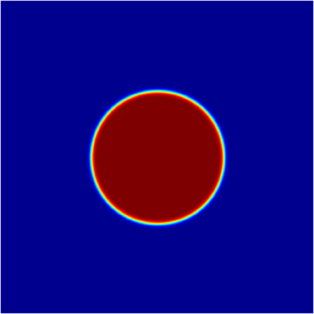}~
\includegraphics[width=0.19\textwidth]{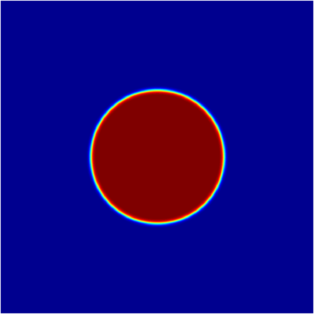}~
\caption{Minimizers of $F_\varepsilon$ for $c=1/7$ and $\varepsilon=1/100,1/150,1/200,1/250,1/300$. The corresponding cost values are: $1.3089,1.3216,1.3276,1.3311,1.3398$}
\label{circle-dec-eps} 

\end{figure}

In the following paragraphs, we take as a toy problem the isoperimetric problem. The third property stated in Proposition \ref{properties-gconv} justifies the following numerical approach. In order to approach the set which minimizes the perimeter at fixed volume, we find minimizers $m_\varepsilon$ of $F_\varepsilon$ for $\varepsilon$ smaller and smaller. We expect that the minimizers $m_\varepsilon$ approach the minimizer of $F$. We consider a straightforward finite differences discretization to compute $F_\varepsilon$ on a fixed grid $N \times N$ in the unit square $[0,1]^2$. The procedure is as follows:
\begin{itemize}
\item Fix an initial $\varepsilon_0$ and a random initial condition, and then compute the numerical minimizer of $F_{\varepsilon_0}$;
\item Decrease $\varepsilon$ and find the numerical minimizer of $F_\varepsilon$ starting from the previous minimizer.
\item Repeat until $\varepsilon$ is small enough.
\end{itemize} 
This simplistic approach has one drawback: the choice of $\varepsilon_0$ cannot be made independent of the grid step. The $\varepsilon$ parameter governs the width of the interface between $0$ and $1$ for the minimizer of $F_\varepsilon$. If $\varepsilon$ is less than $1/N$ then the gradient term in $F_\varepsilon$ contains meaningless information, since the width of the interface is smaller than the width of the grid. To fix this issue, we start with $\varepsilon_0 \in [1/N,4/N]$ and whenever we decrease $\varepsilon$ we refine the grid and interpolate the initial condition on this new grid. We present the numerical results obtained using this procedure in the case $c=1/7$. In this case, we know that in two dimensions, the solution of the isoperimetric problem is a disk, and the corresponding perimeter to a disk of area $1/7$ is $2 \sqrt{\pi/7} = 1.3398$. Results can be seen in Figure \ref{circle-dec-eps}. It is interesting to note that as $\varepsilon$ becomes smaller and smaller, the  minimal values of the functionals $F_\varepsilon$ converge towards the minimal value of $F$, as expected.

\begin{figure}
\includegraphics[width=0.19\textwidth]{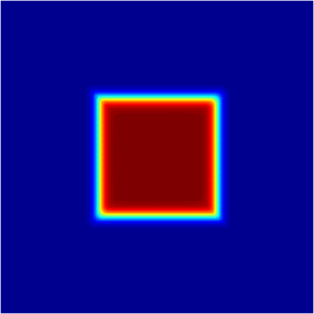}~
\includegraphics[width=0.19\textwidth]{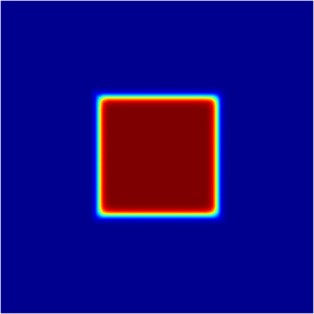}~
\includegraphics[width=0.19\textwidth]{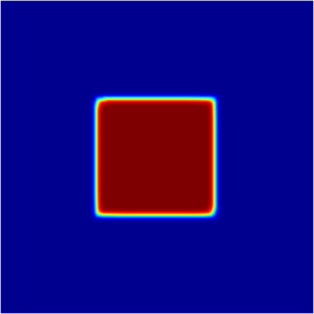}~
\includegraphics[width=0.19\textwidth]{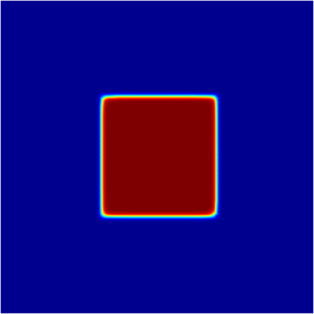}~
\includegraphics[width=0.19\textwidth]{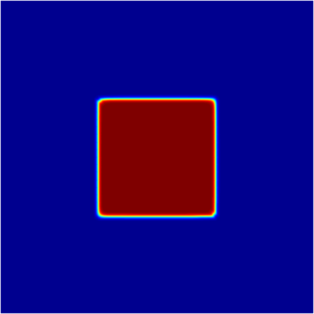}~
\caption{Minimizers of $G_\varepsilon$ for $c=1/7$ and $\varepsilon=1/100,1/150,1/200,1/250,1/300$. The corresponding cost values are: $1.4851,1.4914,1.4979,1.5031,1.5049$}
\label{aniso-dec} 
\end{figure}

We can consider the same problem in an anisotropic setting. If we consider a set $\Omega \subset \Bbb{R}^d$ with $C^1$ boundary, then its perimeter is equal to
\[ \Per(\Omega)=\int_{\partial \Omega}1 d\mathcal{H}^{n-1} =\int_{\partial \Omega} \|\vec n(x)\| d\mathcal{H}^{n-1} \]
where $\vec n(x)$ denotes the unit outer normal vector corresponding to $x \in \partial \Omega$. Thus, the perimeter treats all directions in the same way and no direction has an advantage over the others. Things change if we pick another norm $\varphi$ on $\Bbb{R}^d$, different from the euclidean one. We can define the anisotropic perimeter associated to a norm $\varphi$ by
\[ \Per_\varphi(\Omega)=\int_{\partial \Omega} \varphi(\vec n) .\]
It is possible to prove a variant of the Modica-Mortola theorem in the anisotropic case. Proofs of this result can be found in \cite{braides2},\cite{braides}. A local variant of this result, where the norm $\varphi$ can also depend on the position of the point can be found in \cite{ambrosio-aniso}. In Section \ref{main-results} we provide a different direct proof of the result.
\begin{thm} Let $D$ be a bounded open set and let $W: \Bbb{R} \to [0,\infty)$ be a continuous function such that $W(z)=0$ if and only if $z \in \{0,1\}$. Consider $\varphi$ a norm on $\Bbb{R}^d$. Denote $c =2\int_0^1 \sqrt{W(s)}ds$. We define $G_\varepsilon,G : L^1(D) \to [0,+\infty]$ by
\[ G_\varepsilon(u) = \begin{cases}\varepsilon \int_D\varphi(\nabla u)^2 +\frac{1}{\varepsilon} \int_D W(u) & u \in H^1(D) \cap X\\ +\infty & \text{otherwise} \end{cases} \]
and
\[ G(u)             = \begin{cases} c\Per_\varphi(u^{-1}(1)) & u \in BV(D; \{0,1\})\cap X \\ +\infty & \text{otherwise}\end{cases}\]
then
\[ G_\varepsilon \stackrel{\Gamma}{\longrightarrow} G\]
in the $L^1(D)$ topology.
\label{modica-mortola-aniso}
\end{thm}
We repeat the same experiment as in the isotropic case. Pick $\varphi(x) = |x_1|+|x_2|$, a norm which favorizes the vertical and horizontal directions. Then the shape which minimizes $\Per_\varphi(\Omega)$ with area constraint, the so-called Wulff shape associated to $\varphi$, is a square. When $c=1/7$ the optimal value is $4/\sqrt{7}=1.5118$. In Figure \ref{aniso-dec} we present the optimizers of $G_\varepsilon$ for decreasing values of $\varepsilon$ and we observe the same convergence behavior. 

The next step is to consider partitioning problems. 
One such method was developed by \'E. Oudet in \cite{oudet} in the case of partitions minimizing the sum of perimeters of the cells in two and three dimensions. The author uses a generalization of the Modica-Mortola theorem to the case of partitions. The partition condition, in this functional case, is realized by imposing that the density functions $u_1,u_2,...,u_n$, corresponding to the cells of the partition, satisfy the relation $u_1+u_2+...+u_n=1$. Note that this last condition is not too difficult to implement from a numerical point of view. With this framework the author was able to recover the result of Hales in the plane. In three dimensions the numerical optimizer was close to the Weaire-Phelan structure.

In the following we provide an extension of this numerical framework in the anisotropic case. First we provide a $\Gamma$-convergence result which generalizes Theorem \ref{modica-mortola-aniso} to the partition case. We underline the fact that the $\Gamma$-convergence is not stable for the sum, so the result is not trivial. In fact, the (LI) property in the definition of the $\Gamma$-convergence comes at once from the one phase case, while the (LS) property requires a bit of work. In order to construct a recovery sequence, we use an approximation result proven by Baldo \cite{baldo}, which states that we can approximate well enough every admissible partition by a polygonal partition.

In the end, we present some numerical computations, for different anisotropy choices, and we find the expected behaviour: partition cells tend to have their boundaries aligned with the favorized directions. Although the theoretical framework is restricted to the case where $\varphi$ is a norm, and thus, is convex, we observe numerically that non-convex anisotropies also produce the expected results and the  rate of convergence is much higher in some cases. 

\section{Main Results}
\label{main-results}

We consider the following definition of the generalized perimeter, valid for every measurable set $D\subset \Bbb{R}^d$.
\[ \Per(\Omega,D) = \sup\{ \int_\Omega \di g dx : g \in C_0^\infty(D;\Bbb{R}^d),\|g\|_\infty \leq 1\} \]
This definition agrees with the classical one in the case $\Omega$ has a certain regularity (polyhedra, piecewise $C^1$, etc). Given a norm $\varphi$ on $\Bbb{R}^d$ we can extend the above variational characterization to the anisotropic perimeter associated to $\varphi$ by
\[ \Per_\varphi(\Omega,D) = \sup\{ \int_\Omega \di g dx : g \in C_0^\infty(D;\Bbb{R}^d),\varphi(g) \leq 1\}. \]
We make the assumption that $\varphi$ is comparable with the Euclidean norm, i.e. there exist constants $c,C>0$ such that $c|x|\leq \varphi(x)\leq C|x|$. Then if a set $E$ satisfies $\Per_\varphi(E,\Omega)<\infty$ then $\chi_E \in BV(\Omega)$, the space of functions of bounded variation on $\Omega$. 

Furthermore, we can choose norms which depend on the position of the considered point: $\varphi : D \times \Bbb{R}^N$ which are lower semicontinuous, positively $1$-homogeneous and convex in the second variable. In addition, we assume the existence of $0<m\leq M$ such that $m |\xi|\leq \varphi(x,\xi) \leq M|\xi|$ for every $(x,\xi) \in D \times \Bbb{R}^N$. Then a local anisotropic perimeter can be defined as follows
\[ \Per_\varphi(\Omega,D) = \sup\{ \int_\Omega \di g dx : g \in C_0^\infty(D;\Bbb{R}^d),\varphi(x,g(x)) \leq 1\}. \]

The purpose of the following paragraphs is to approximate by $\Gamma$-convergence the sum of the anisotropic perimeters of a partition of a bounded, open set $D$ into $n$ parts of equal volumes. We want to be able to have a result which is also valid for local anisotropies, where the norm $\varphi$, which determines the anisotropy may also depend on the point $x$. The $\Gamma$-convergence result is divided in two parts, corresponding to the two properties in its definition. The (LI) property can be deduced by studying the one phase case. The (LS) property needs some work in order to construct a suitable recovery sequence.

The double-well potential $W$ is stated in a general form in the theorem, but we will assume that it has additional properties. In practice we use $W(s)=s^2(1-s)^2$, but we are only interested of the form of $W$ in a neighborhood of $[0,1]$. Therefore, we assume that $W$ is bounded (by truncating it at a large level, if necessary). In order to simplify the construction of the recovery sequence, we assume that the graph of $W$ is symmetric with respect to the line $x = 1/2$. The theorem stated below is a particular case of the one studied in \cite{ambrosio-aniso}. We give a slightly different proof, and adapt it to the case of partitions.

\begin{thm}
Let $D$ be an open, bounded domain in $\Bbb{R}^N$, and $f : D \times \Bbb{R}^N \to [0,\infty]$ be a lower semicontinuous function, positively $1$-homogeneous and convex in the second variable, which satisfies $m|\xi| \leq \varphi(x,\xi)\leq M|\xi|$ for every $(x,\xi) \in D \times \Bbb{R}^N$, with $0<m\leq M$. We consider $W : \Bbb{R} \to [0,\infty)$ such that $W(0)=W(1) = 0$ and $W(x)>0$ for $x \notin \{0,1\}$. Define $F_\varepsilon,F : L^1(D) \to [0,\infty]$ as follows:
\[ F_\varepsilon(u) = \begin{cases}
\ds\varepsilon \int_D \varphi(x,\nabla u(x))^2dx+\frac{1}{\varepsilon} \int_D W(u(x))dx & \text{ if } \ds u \in H^1(D), \ \int_D u = c\\
+\infty & \text{ otherwise}
\end{cases}
\]
\[ F(u) = \begin{cases} 
c\ds \int_{S(u)} \varphi(x,\nu_u)d\mathcal{H}^{N-1} & \text{ if } u \in BV(D,\{0,1\}),\ \ds \int_D u = c\\
+\infty & \text{ otherwise}
\end{cases}, 
\]
where $c = 2\int_0^1 W(s)^{1/2}ds$ and $S(u)$ is the jump set of $u$.

Then for every $u \in L^1(D)$ and every $(u_\varepsilon) \in L^1(D)$ such that $(u_\varepsilon) \to u$ in $L^1(D)$ we have
\[
\liminf_{\varepsilon}
F_\varepsilon(u_\varepsilon) \geq F(u).
\]
\label{single}
\end{thm}

\emph{Proof:} This result follows naturally from the following remarks and from a variant of Reshetnyak's semicontinuity theorem. Consider the function $\phi(t) = 2\int_0^t W(s)^{1/2}ds$, which is Lipschitz continuous, in view of the fact that we assume that $W$ is bounded above. In the following we show that $F(u) = \int_D \varphi(x,D(\phi\circ u))$, where we use the notation 
\[ \int_D \varphi(x,\mu) = \int_\Omega \varphi\left(x,\frac{d\mu}{d|\mu|}\right) d|\mu|,\]
for every Radon measure $\mu \in \mathcal{M}(D,\Bbb{R}^N)$. First note that if $u \in BV(D,\{0,1\})$ then using the definition of the variation of a $BV$ function we can see that $D(\phi \circ u) = \phi(1) Du$. Moreover, if we have a function $u \in BV(D)$ whose image contains only two real values, then the absolutely continuous part and the Cantor part of $Du$ are zero, while the jump part is
\[ D^j u(B) = \int_{B \cap S(u)} (u^+-u^-)\nu_u d\mathcal{H}^{N-1}\]
where $\nu_u$ is the normal to the jump set $S(u)$ defined by $Du = \nu_u|Du|$. In this case, where $u\in \{0,1\}$ a.e. we also have $Du = d\mathcal{H}^{N-1}\mrestr S(u)$. For details see \cite{braides2}. Having these in mind and using the fact that $\varphi$ is homogeneous of degree one in the second variable, we obtain
\begin{align*}
&c\int_{S(u)} \varphi(x,\nu_u)d\mathcal{H}^{N-1} = c \int_\Omega \varphi\left(x,\frac{dDu}{d|Du|}\right)d\mathcal{H}^{N-1}\mrestr S(u) \\
=& \phi(1) \int_\Omega \varphi \left(x,\frac{dDu}{d|Du|}\right)d |Du|=\int_\Omega \varphi \left(x,\frac{dD(\phi\circ u)}{d|D (\phi\circ u)|}\right)d |D(\phi\circ u)|\\
&= \int_D \varphi(x,D(\phi\circ u)).
\end{align*}

The following variant of Reshetnyak lower semicontinuity theorem can be found in \cite[Theorem 2.38]{ambrosiofuscopallara}.

\begin{thm}
Let $D$ be an open subset of $\Bbb{R}^N$ and $\mu,\mu_n$ be $R^n$-valued finite Radon measures in $D$. If $\mu_n \to \mu$ weakly* in $D$ then
\[ \int_D f(x,\mu) \leq \liminf_{ n \to \infty} \int_D f(x,\mu_n),\]
for every lower semicontinuous function $f: \Omega \times \Bbb{R}^n \to [0,\infty]$, positively $1$-homogeneous and convex in the second variable.
\label{reshetnyak}
\end{thm}

First, let's note that the integral condition is preserved under $L^1(D)$ convergence, since
\[ \left| \int_D u_\varepsilon - \int_D u \right| \leq  \|u_\varepsilon-u\|_{L^1(D)}.\]  Since $\phi$ is Lipschitz continuous, $u_\varepsilon \to u$ in $L^1(D)$ implies that $\phi \circ u_\varepsilon \to \phi \circ u$ in $L^1(D)$. If we suppose that $ \liminf_{\varepsilon \to 0} F_\varepsilon(u_\varepsilon)<+\infty$ (else there is nothing to prove) then, using the standard inequality $a+b \geq \sqrt{ab}$, we get that $F_\varepsilon(u_\varepsilon) \geq 2 \int_D \varphi(x,D(\phi\circ u_\varepsilon)) \geq 2 \int_D D(\phi\circ u_\varepsilon)$. Therefore, we can assume that $\sup |D(\phi\circ u_\varepsilon)|(D)<+\infty$. According to \cite[Definition 1.41, Remark 1.42]{braides2} we can conclude that $D(\phi\circ u_\varepsilon) \wconv D(\phi \circ u)$ weakly* in $\mathcal{M}(D,\Bbb{R}^N)$ and Theorem \ref{reshetnyak} is applicable:
\begin{align*}
& \liminf_{\varepsilon \to 0}F_\varepsilon(u_\varepsilon) \geq \liminf_{\varepsilon \to 0} 2 \int_D \varphi(x,\nabla u_\varepsilon)W(u_\varepsilon)^{1/2}\\
=&\liminf_{\varepsilon\to 0} \int_D \varphi(x,D(\phi\circ u_\varepsilon))\geq \int_D \varphi(x,D(\phi \circ u))=F(u).
\end{align*}
\hfill $\square$

We are now able to state the $\Gamma$-convergence result concerning the partition case. We use the notation. We assume that the potential $W$ satisfies the following properties:
\begin{itemize}
\item $W$ satisfies the hypotheses of Theorem \ref{single}.
\item $W(0.5-t) = W(0.5+t)$ for every $t \in \Bbb{R}$.
\item $W$ is bounded above.
\end{itemize}
We also assume that $\varphi : D \times \Bbb{R}^N \to [0,\infty)$ satisfies the hypotheses of Theorem \ref{single} and that it is Lipschitz continuous in the first variable. We use the following bold notation to denote vectors of functions: $\bo{u}=(u_i) \in L^1(D)^n$.

In the following, we consider $X \subset L^1(D)^n$ to be the space containing the $n$-uples of function satisfying the partition condition and the area constraints:
\[ X = \{ \bo{u}  \in L^1(D)^n : \int_{D}u_i = \frac{|D|}{n},\ u_1+...+u_n = 1 \text{ in }D\}.\]
We note the fact that the proofs which follow do not change much if instead of the equal areas conditions we put only a fixed area condition on every one of the phases.
\begin{thm}
We consider the functionals $F_\varepsilon,F : (L^1(D))^n \to [0,\infty]$, defined by
\[ F_\varepsilon(\bo{u})= \begin{cases} \ds 
\sum_{i=1}^n \left( \varepsilon \int_D \varphi(x,\nabla u_i)^2 +\frac{1}{\varepsilon} \int_D W(u_i) \right) & \text{ if } \bo{u} \in (H^1(D))^n \cap X \\
 +\infty & \text{ otherwise}
\end{cases}\]
\[ F(\bo{u})=\begin{cases}
\sum_{i=1}^n c\ds \int_{S(u_i)} \varphi(x,\nu_{u_i}) & \text{ if } \bo{u} \in (BV(D,\{0,1\})^n \cap X\\
+\infty & \text{ otherwise }
\end{cases}
\]
Then $F_\varepsilon \gconv F$ in the $(L^1(D))^n$ topology.
\label{aniso-partitions}
\end{thm}

\emph{Proof:} The (LI) part of this result follows at once from Theorem \ref{single}.

For the (LS) part we need to be able to construct a recovery sequence for every $\bo u \in L^1(D)$ such that $F(\bo u)<+\infty$. In order to do this, we reduce the problem to subset $\mathcal D \subset \{F<+\infty\}$ which is dense and has some good regularity properties. This is a classical procedure described in Proposition \ref{reduction-regular} and \cite{braides2}. One such suitable dense class is provided by Baldo in \cite{baldo} and consists of functions $\bo u \in BV(D,\{0,1\})^n \cap X$ which represent partitions of $D$ into polygonal domains. 

The result of Baldo says that for every $\bo u \in (BV(D,\{0,1\})^n \cap X$ there exists a sequence $\bo u_n \in (BV(D,\{0,1\})^n \cap X$ such that $\bo u_n \to \bo u$ in $(L^1(D))^n$, each component of $\bo u_n$ represents a set of finite perimeter, $D\bo u_n^i \wconv D\bo u^i$ weakly* in $\mathcal{M}(D,\Bbb{R}^N)$ and $|D \bo u_n^i|(D) \to |D\bo u_n|(D)$ (the corresponding perimeters converge). The Reshetnyak continuity theorem found in \cite[Theorem 2.39]{ambrosiofuscopallara} assures us that $F(\bo u_n) \to F(\bo u)$. Thus Proposition \ref{reduction-regular} allows us to restrict our attention to functions $\bo u$ which represent partitions of $D$ into polygonal domains of equal areas.

We consider the optimal profile problem 
\[ c = \min\left\{ \int_R (W(v)+|v'|^2) dt : v(-\infty) = 0 , v(+\infty)=1 \right\}\]
and the related problem
\begin{equation} zc = \min\left\{\int_R (W(v)+z^2|v'|^2) dt : v(-\infty) = 0 , v(+\infty)=1 \right\}
\label{optprof2}
\end{equation}
Note that the solution of \eqref{optprof2} satisfies the differential equation $v' = \sqrt{W(v)}/z$ and for symmetry reasons, we impose the initial condition $v(0)=1/2$. Note that $v$ is strictly increasing, and $v(t) \geq 1/2$ for $t \geq 0$. It is not difficult to see that $c = \ds 2 \int_0^1 \sqrt{W(s)}ds$.

Take $v$ a solution to problem \eqref{optprof2}.  We modify $v$ such that it goes from $0$ to $1$ on a finite length interval in the following way (inspired from \cite{braides2}):
\[ v^\eta = \min\{\max\{0,(1+2 \eta)v -\eta\},1\}.\] 
We have 
\[ c^\eta = \int_\Bbb{R} (W(v^\eta)+|(v^\eta)'|^2) \to c \text{ as } \eta \to 0.\]

\begin{figure}
\centering
\includegraphics[width = 0.4\textwidth]{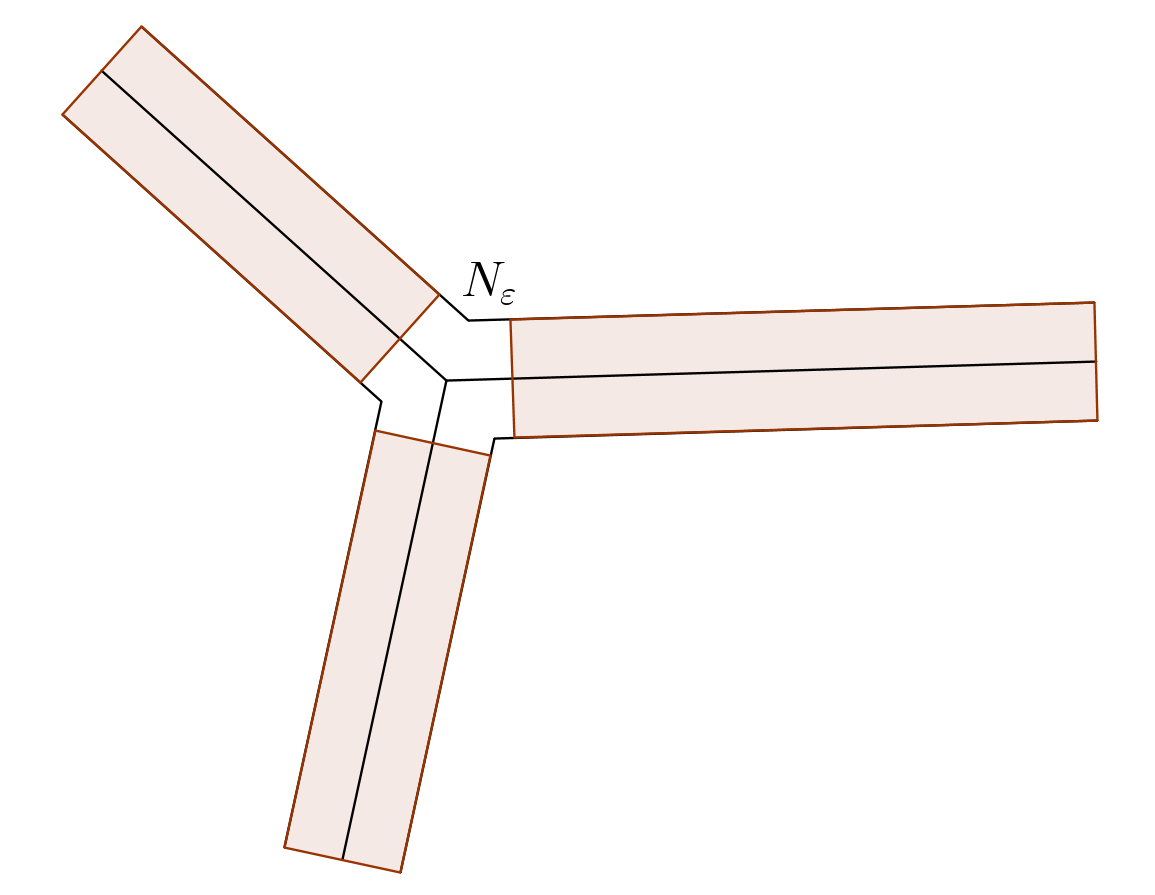}
\caption{Example of a part of $N_\varepsilon$}
\label{rectangles}
\end{figure}

We denote $(\Omega_i)_{i=1}^n$ the polygonal partition determined by $\bo u$. We denote by $N_\varepsilon$ the set 
of points which are close to triple (or multiple) points of the partition $(\Omega_i)$, such that 
\[ \{ x \in D : d(x,\Omega_i)<\varepsilon\} \setminus N_\varepsilon,\]
is a union of rectangles. An example is given in Figure \ref{rectangles}.

In the following, we denote by $v_{\vec{n}}$ the optimal profile with $z = \varphi(\vec{n})$. We use the signed distance $d_E(x) = d(x,D \setminus E)-d(x,E)$ and define $u_\varepsilon^i$ on $D \setminus N_\varepsilon$ by
\[ u_\varepsilon^i (x) = \begin{cases}
v_{\nabla{d_{\Omega_i}(x)}}^\eta (\frac{d_{\Omega_i}(x)}{\varepsilon}) & \text{ if } |d_{\Omega_i}(x)|\leq T\varepsilon \\
0 & \text{ otherwise in }  D \setminus \Omega_i \\
1 & \text{ otherwise in }  \Omega_i
\end{cases}
\]
where $T$ is great enough such that the support of $(v^\eta)'$ is contained in $[-T,T]$. Until here, $u_\varepsilon$ is a Lipschitz continous function with values in $[0,1]$ and a Lipschitz constant of order $1/\varepsilon$. We extend each $u_\varepsilon^i$ to the whole $D$ with the same Lipschitz constant (this is possible by Kirszbraun's theorem, see \cite{gariepy}).

In order that $\bo u_\varepsilon \in X$ we must treat the measure and the sum constraints. We deal with the sum constraint first. We have three types of points:
\begin{itemize}
\item $|d_{\Omega_i}(x)| > T\varepsilon$ for all $i$. Here the sum constraint is clear, since one component takes value $1$ and the rest $0$.
\item There exist precisely $2$ indexes $i,j$ such that $|d_{\Omega_i}(x)|,|d_{\Omega_j}(x)|\leq T\varepsilon$. Here the symmetry of the optimal profile assures us that the $u_\varepsilon^i(x)+u_\varepsilon^j(x)=1$, while the other components take the value $0$.
\item The points in $N_\varepsilon$.
\end{itemize}
We see that the only problems that can occur take place in $N_\varepsilon$. Here, we replace $u_\varepsilon^i$ by $ u_\varepsilon^i/(\sum_{j=1}^n u_\varepsilon^j)$. This operation is well defined, since each $u_\varepsilon^i$ is greater than $1/2$ on $\Omega_i$; thus their sum is always greater than $1/2$. Furthermore, doing this change still leaves the gradient of $u_\varepsilon^i$ of the form $O(1/\varepsilon)$.

In the following we ommit the substript from $v_{\nabla d_{\Omega_i}(x)}$, and we may do so without loss of generality, since the inequalities described below do not use this dependence until the last few inequalities. Because of the fact that $u_\varepsilon^i$ varies only in the direction of the normal to $\Omega_i$ on $D\setminus N_\varepsilon$, we find that $\nabla u_\varepsilon^i(x)/|\nabla u_\varepsilon^i(x)|$ is a unit normal to $\Omega_i$.

 The integral constraints can be imposed in one of the following ways:
 
 \begin{itemize}
 \item by slightly moving the initial boundaries of $(\Omega_i)$ and then performing the algorithm described above.
 \item by performing the procedure described in \cite{modica}. We modify each phase in a ball of fixed, small enough size, which depends on $\varepsilon$ in order to fix the volume constraints. In the end we note that these perturbations vanish in the limit. 
 \end{itemize}

We split the (LS) estimate in two parts, one on $N_\varepsilon$ and one on $D \setminus N_\varepsilon$.
\begin{align*}
& \int_{N_\varepsilon}\left( \varepsilon \varphi(x,\nabla u_\varepsilon^i)^2+\frac{1}{\varepsilon}W(u_\varepsilon^i)\right) \\
\leq & \frac{|N_\varepsilon| \max_{[0,1]}W}{\varepsilon}+\frac{|N_\varepsilon| \varepsilon \sup_\Bbb{R} |(v^\eta)'|^2 \sup_{\|\vec{n}\|=1}\varphi(x,\vec{n})^2}{\varepsilon} = O(\varepsilon),
\end{align*}
since $|N_\varepsilon|=O(\varepsilon^2)$. This proves that the part corresponding to $N_\varepsilon$ is negligible int the (LS) estimate.

We continue our estimate on $D \setminus N_\varepsilon$: 
\begin{align*}
&\int_{D \setminus N_\varepsilon} \left( \varepsilon \varphi(x,\nabla u_\varepsilon^i)^2+\frac{1}{\varepsilon}W(u_\varepsilon^i) \right) \\
=&\int_{D \setminus N_\varepsilon} \left( \varepsilon \varphi^2(x,\nabla u_\varepsilon/|\nabla u_\varepsilon|)|\nabla u_\varepsilon|^2 
+\frac{1}{\varepsilon} W(v^\eta(d_{\Omega_i(x)}(x))/\varepsilon)\right)\\
= & \int_{-T\varepsilon}^{T\varepsilon} \int_{\{d(x)=t\}\setminus N_\varepsilon}
\left( \varepsilon \varphi^2(x,\nu_{\bo u^i}) \frac{|(v^\eta)'(t/\varepsilon)|^2}{\varepsilon^2} +\frac{1}{\varepsilon} W(v^\eta(d_{\Omega_i(x)}(x))/\varepsilon) \right)d \mathcal{H}^{N-1}(x) dt\\
= & \int_{S(\bo u_\varepsilon^i) \setminus N_\varepsilon}  \int_{-T\varepsilon}^{T\varepsilon} \left(\frac{1}{\varepsilon} W(v^\eta(t/\varepsilon))+\frac{1}{\varepsilon}\varphi^2(x,\nu_{\bo u_i}(x))|(v^\eta)'(t/\varepsilon)|^2 \right)dt d\mathcal{H}^{N-1}(x)+O(\varepsilon)\\
= & \int_{S(\bo u_\varepsilon^i) \setminus N_\varepsilon}  \int_{-T}^{T} \left( W(v^\eta(t))+\varphi^2(x,\nu_{\bo u_i}(x))|(v^\eta)'(t)|^2 \right)dt d\mathcal{H}^{N-1}(x)+O(\varepsilon)\\
\leq & c^\eta \int_{S(\bo u^i)} \varphi(x,\nu_{\bo u_i}) d\mathcal{H}^{N-1}+O(\varepsilon). 
\end{align*}
We have used the co-area formula. The fact that $\varphi$ is Lipschitz continuous in the first variable allows us to write estimates of the form $\varphi(y,\xi) \leq \varphi(x,\xi)+L|x-y|$, and this is why we have an $O(\varepsilon)$ term after we change the order of integration. The (LS) property comes from summing the estimates obtained for every $(\bo u_\varepsilon^i)$.

\section{Numerical Computations}

\subsection{Partitions minimizing an anisotropic perimeter}

One of the main properties of the $\Gamma$-convergence is the fact that if $F_\varepsilon \gconv F$ then any limit point of a sequence $(x_\varepsilon)$ of minimizers of $F_\varepsilon$ is a minimizer for $F$. Based on this property, we assume that minimizing $F_\varepsilon$ for $\varepsilon$ small enough will get us close to a minimizer of $F$.

We want to approximate numerically the partitions which minimize the sum of their anisotropic perimeters, with respect to some anisotropy $\varphi$. In order to do this, we search numerically for minimizers of 
\begin{equation}
 F_\varepsilon(\bo{u})= 
\sum_{i=1}^n \left( \varepsilon \int_D \varphi(x,\nabla u_i)^2 +\frac{1}{\varepsilon} \int_D W(u_i) \right)
\label{relaxation}
\end{equation}

Using the fact that $\varphi(x,\xi) \geq c|\xi|$ for a constant $c>0$, we deduce that if $\bo{u_n}$ is a minimizing sequence for $F_\varepsilon$ then $(\nabla \bo{u}_n^i)$ is bounded in $L^2(D)$. Truncating $(\bo{u}_n)$ between $0$ and $1$ decreases $F_\varepsilon(\bo{u}_n)$, so $(\bo{u}_n)$ is also bounded in $L^2(D)^n$. Thus $(\bo{u}_n)$ is bounded in $H^1(D)^n$, which means that it has a subsequence which converges weakly $H^1$ to $\bo u$. The convexity of $\varphi$ and the Fatou Lemma imply that
\[ \liminf_{n \to \infty}F_\varepsilon (\bo u_n) \geq F(\bo u),\]
which means that \eqref{relaxation} has a minimizer in $H^1(D)^n$. The lack of convexity of the potential $W$ does not allow us to conclude that the minimizer is unique. In fact, domain symmetry and permutations of phases always lead to multiple optimizers.

We can devise an algorithm to approximate numerically such a minimizer. We discretize the unit square $D=[0,1]^2$ using a finite differences grid, and use quadrature formulas to compute the integrals in the expression of $F_\varepsilon$. The choice of $\varepsilon$ is important in order to have meaningful results. Morally, $\varepsilon$ dictates the width of the interface between the sets $\{u_i = 0\}$ and $\{u_i = 1\}$, and it cannot be lower than the width of the discretization grid. Satisfactory results have been obtained for $\varepsilon \in [\frac{1}{N},\frac{4}{N}]$. Note that if $\varepsilon$ is large then the diffusion interface is bigger, and therefore the shapes can move more freely in order to find their optimal position. Forcing $\varepsilon$ small in the beginning may lead to a local minimum.
In order to diminish the size of the interface, we can iterate the optimization algorithm by decreasing $\varepsilon$.

We observe that the behavior of the algorithm depends heavily on the choice of $\varphi$. We have many options to choose the anisotropy $\varphi$: 
\begin{itemize}
\item $\varphi(x) = |x_1|+|x_2|$ - horizontal and vertical directions;
\item $\varphi(x) = \left(|x_1|^p+|x_2|^p\right)^{1/p}$
\item $\varphi(x) = |ax_1+bx_2|+|cx_1+dx_2|$ - variable directions corresponding to $a,b$.
\item $\varphi(x) = (ax_1^2+bx_2^2)^{1/2}$ with $a>b$: favorize one of the directions corresponding to coordinate axes.
\end{itemize}
We present below some numerical results we obtained using various norms and parameters.

\begin{figure}
\centering
\includegraphics[width = 0.3\textwidth]{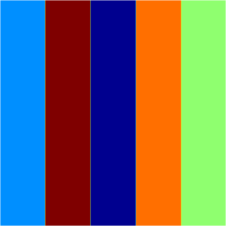} \includegraphics[width = 0.3\textwidth]{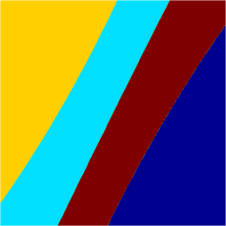}
\includegraphics[width = 0.3\textwidth]{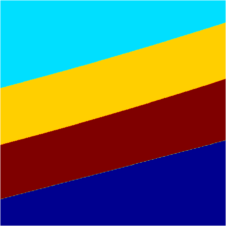}
\caption{Examples of optimal partitions with one favorized direction}
\label{aniso-1dir}
\end{figure}

The first example we study is the case where we have one favorized direction. Favorizing one direction parallel to the coordinate axis is not hard. It is enough to use a weighted norm like $\varphi(x) = \sqrt{x_1^2+100x_2^2}$ to favorize the vertical direction. Indeed, looking at the term $\ds \int_D \varphi(\nabla u)$ we see that if the gradient $\nabla u$ has a second component which is large, then the quantity $\varphi(\nabla u)$ is large. Thus, in order to minimize our functional, the gradient of $u$ should be close to zero in the second component. Thus $u$ is close to a constant on each vertical line, and all boundaries will be vertical at the optimum. In order to favorize a general direction, one could use a rotation of the coordinate axis included in the norm. A few examples of optimal partitions with one favorized direction can be seen in Figure \ref{aniso-1dir}.

The next interesting situation is the case of two favorized directions. Since we work on rectangular domains, it is natural to consider vertical and horizontal favorized directions. This can be achieved using the $\ell^1$ norm $\varphi(x) = |x_1|+|x_2|$. Another way of favorizing these two direction is presented below.

\begin{figure}
\centering
\includegraphics[width = 0.18\textwidth]{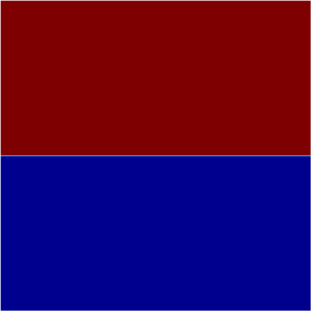}
 \includegraphics[width = 0.18\textwidth]{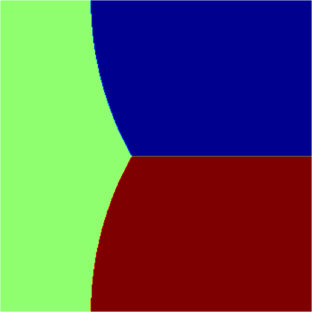} 
\includegraphics[width = 0.18\textwidth]{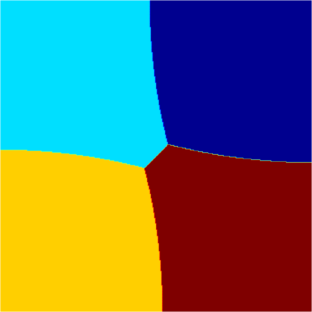}
\includegraphics[width = 0.18\textwidth]{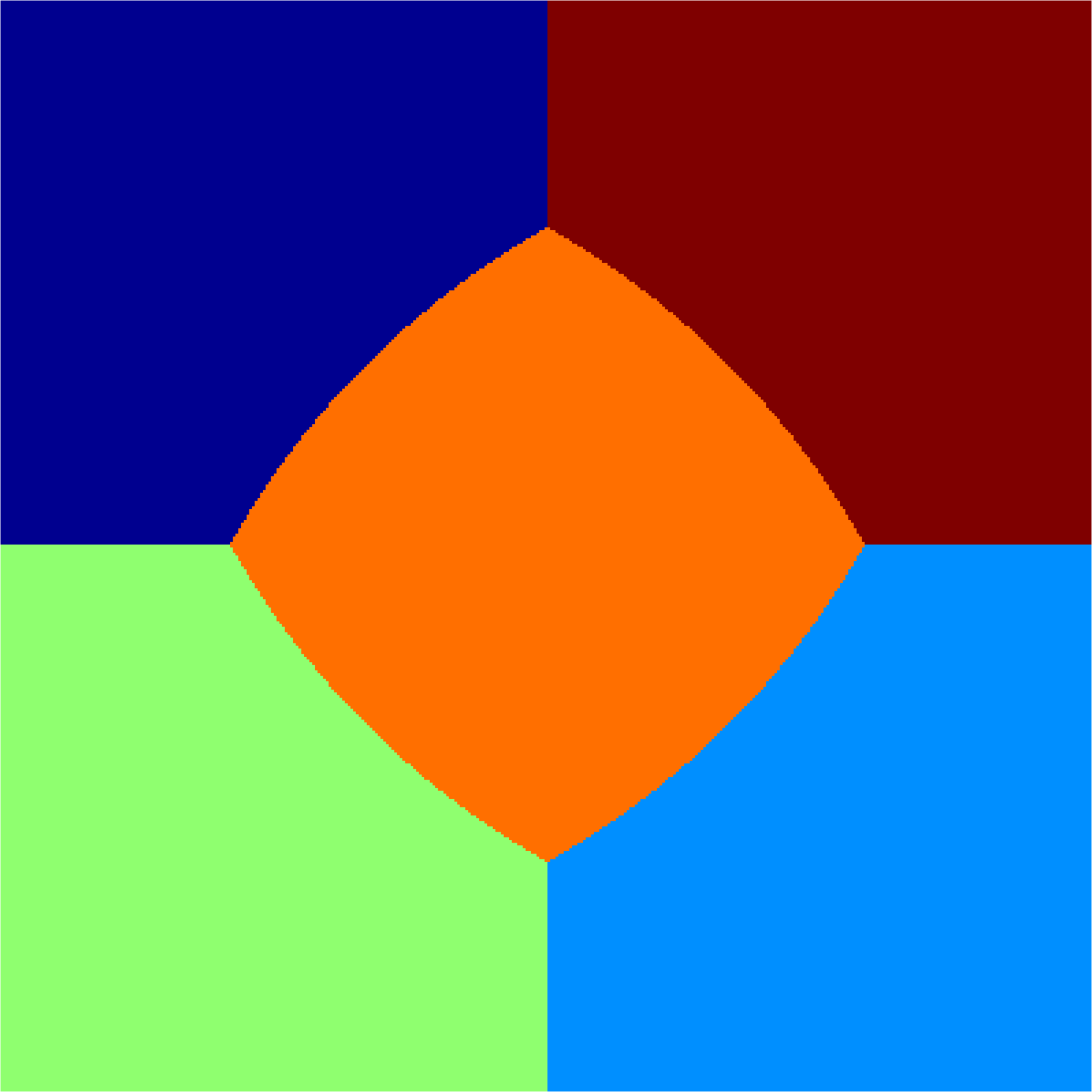}
 \includegraphics[width = 0.18\textwidth]{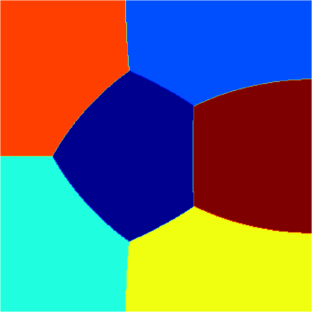}
 
 \vspace{0.1cm}
 
\includegraphics[width = 0.18\textwidth]{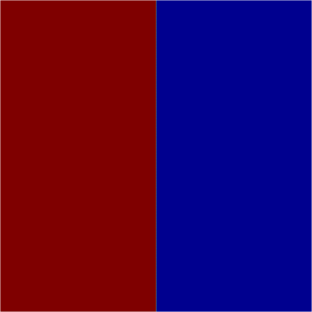}
 \includegraphics[width = 0.18\textwidth]{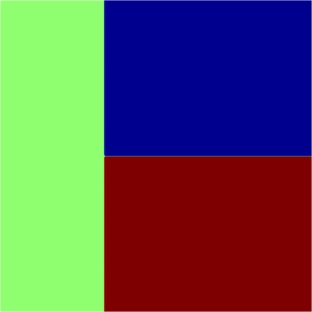} 
\includegraphics[width = 0.18\textwidth]{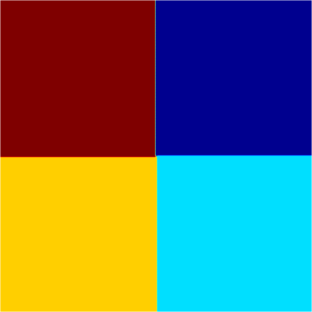}
\includegraphics[width = 0.18\textwidth]{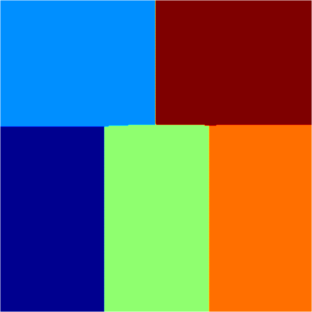}
 \includegraphics[width = 0.18\textwidth]{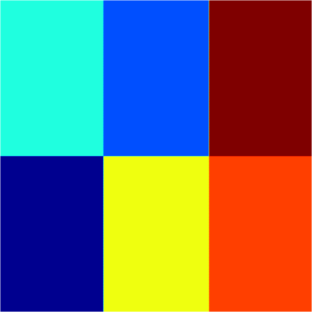}
 
 \vspace{0.2cm}
 
\includegraphics[width = 0.18\textwidth]{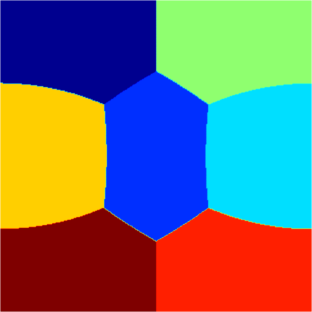}
\includegraphics[width = 0.18\textwidth]{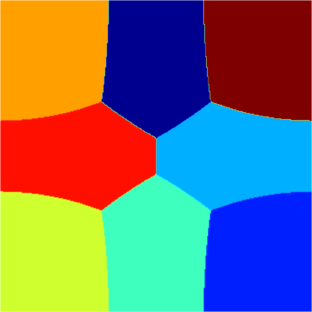} \includegraphics[width = 0.18\textwidth]{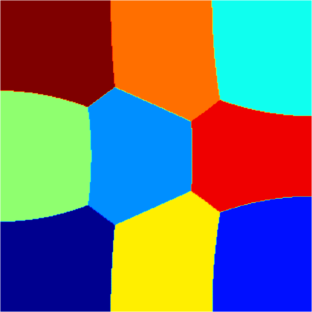}
\includegraphics[width = 0.18\textwidth]{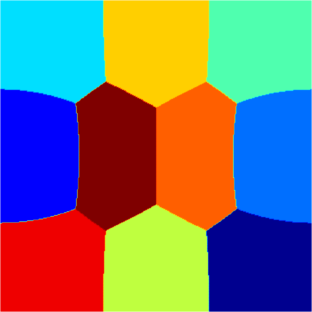}

\vspace{0.1cm}

\includegraphics[width = 0.18\textwidth]{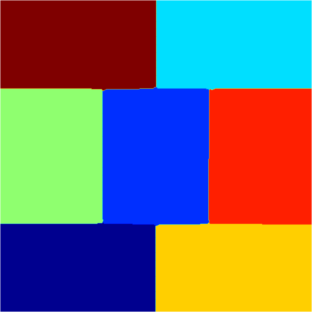}
\includegraphics[width = 0.18\textwidth]{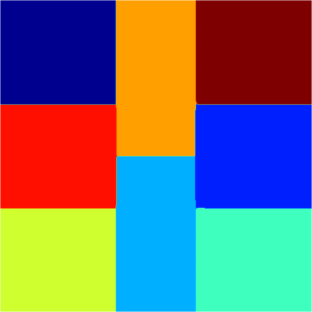} \includegraphics[width = 0.18\textwidth]{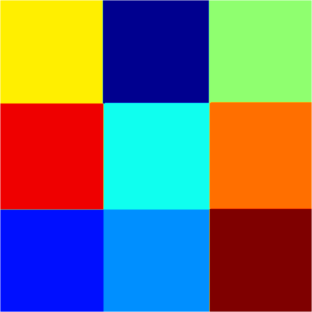}
\includegraphics[width = 0.18\textwidth]{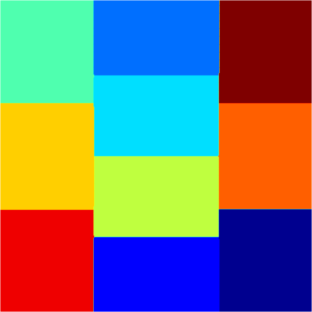}
\caption{Optimal partitions obtained for $N\in [2,10]$. The isotropic case (up) and the anisotropic case corresponding to $\varphi(x) = |x_1|+|x_2|$ (down)}

\label{aniso-comp}
\end{figure}

One natural way to favorize a direction corresponding to a coordinate axis is to use a norm of the form 
\[ \varphi(x) = \sqrt{ax_1^2+bx_2^2},\]
with $a>b$. In order to favorize two directions we can think of using something of the form
\[ \varphi(x) = \sqrt[4]{(100x_1^2+x_2^2)(x_1^2+100x_2^2)}.\]
The problem with the above choice of $\varphi$ is the lack of convexity, which goes out of the $\Gamma$-con\-ver\-gence framework of the theoretical result. Nevertheless, we observe that despite this non-convexity issue we obtain the same results as in the case of the $\ell^1$ norm. Moreover, the convergence is accelerated in the case of the non-convex $\varphi$. We present in Figure \ref{aniso-comp} the partitions of the unit square corresponding to the classic perimeter and the ones obtained favorizing horizontal and vertical directions. Since the results we obtained are all partitions of the square in rectangles of equal areas, we may ask if these rectangle configurations are optimal. The answer is yes, and the problem of partitioning a square into rectangles of equal areas which minimize their total perimeter has been completely answered in \cite{partition-square}.

\begin{figure}
\centering
\includegraphics[width = 0.23\textwidth]{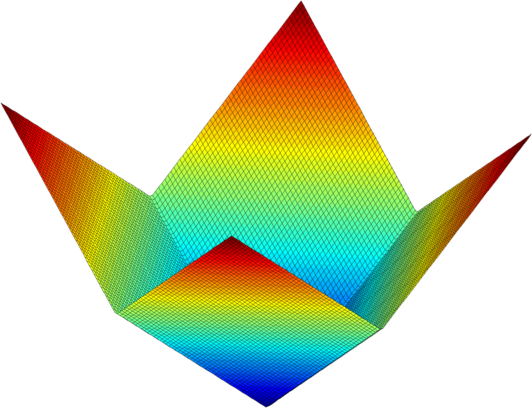}
 \includegraphics[width = 0.23\textwidth]{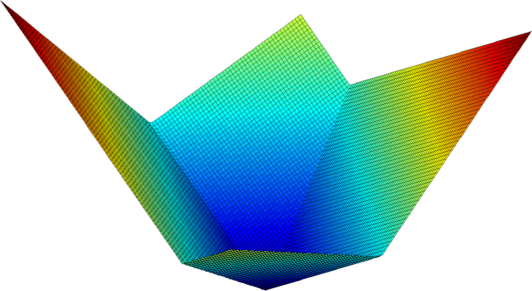} 
\includegraphics[width = 0.23\textwidth]{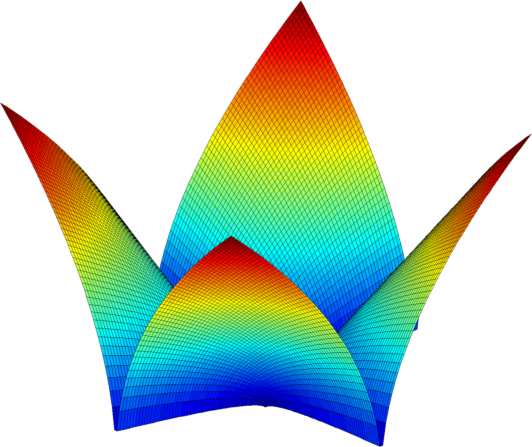}
\includegraphics[width = 0.23\textwidth]{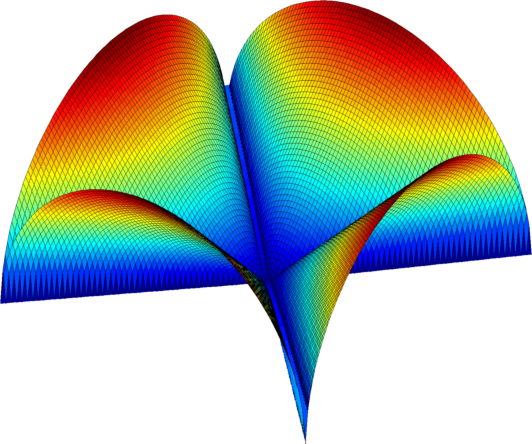}
\caption[Plots of some of the norms we considered, on the unit square.]
    {\tabular[t]{@{}l@{}}
    Plots of some of the norms we considered, on the unit square. In order left to right: \\ 
     1. $\ell^1$ norm, directions $0,\pi/2$. \\
     2. $\ell^p$ norm, $p=1.1$, directions $\pi/6,\pi/2$. \\
     3. Square root of product of two norms, directions $0,\pi/2$.\\
     4. Square root of product of two norms, directions $-\pi/4,\pi/4$.
      \endtabular} 
\label{norm-plots}
\end{figure}

\begin{figure}
\centering
\includegraphics[width = 0.3\textwidth]{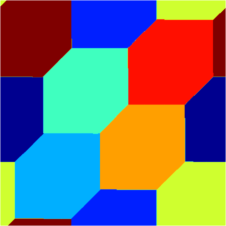} \includegraphics[width = 0.3\textwidth]{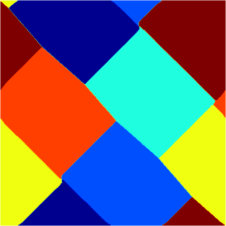}
\includegraphics[width = 0.3\textwidth]{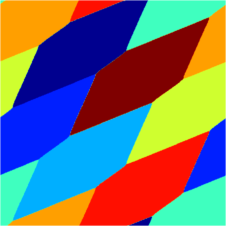}
\caption{Optimal partitions for other anisotropies with two or three favorized directions, under periodicity conditions}
\label{aniso-3dir}
\end{figure}

 As in the case of one favorized direction, we can favorize any desired direction by introducing a suitable rotation in the formulation of the norm. For example, one can favorize the directions corresponding to the two axis bisectors by considering 
\[ \varphi(x) = |x_1+x_2|+|x_1-x_2|.\]
We can continue our study by considering three favorized directions. The choice of the norms is similar, but involving three directions instead of two. As before, we notice a faster convergence when considering non-convex variants of $\varphi$. This behavior could be attributed to the fact that in the non-convex case, the boundaries align immediately to the favorized directions, since along these directions the functional has much lower values. In Figure \ref{norm-plots} you can see some plots of some of functions $\varphi$ we considered, on the unit square. In these picture you can clearly see the favorized directions as the directions along which the lowest values can be found. In the non-convex cases, these directions are more emphasized. Some further computations involving cases where we have three favorized directions can be found in Figure \ref{aniso-3dir}.

\begin{figure}
\centering
\includegraphics[width = 0.3\textwidth]{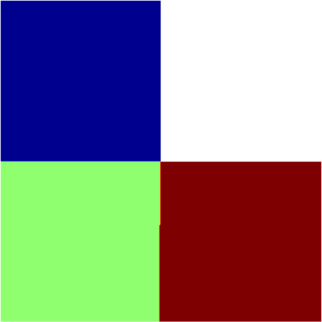} \includegraphics[width = 0.3\textwidth]{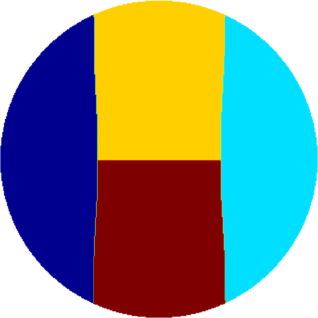}
\includegraphics[width = 0.3\textwidth]{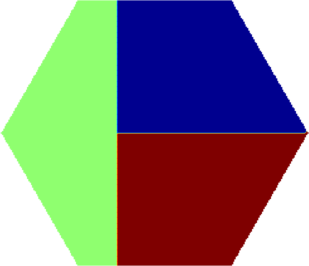}
\caption{Examples of optimal partitions with two favorized directions on general non-rectangular domains}
\label{aniso-general}
\end{figure}

We can use the finite difference framework in the case of non rectangular domain in the following way. We consider the general domain $D$ as a subset of a rectangular region $R$. On this rectangular region a finite differences grid is considered. We apply the same algorithm with the difference that we ignore the grid points which are outside the domain $D$, by assigning them a fixed value zero for the density function and for the gradient of this function. The computation results are not always well behaved near the boundary of $D$, as expected. We present some of the results obtained on general domains in Figure \ref{aniso-general}.

\subsection{Weighted isoperimetric problems}

We note that if the anisotropy functional $\varphi(x,\xi)$ is independent of the direction $\xi$ then the $\Gamma$-convergence result proved in Section \ref{main-results} gives an approximation for the density perimeter $\Per_\nu(\Omega) = \int_{\partial \Omega} \nu(x)d\mathcal{H}^{N-1}$. There are many recent works which treat problems concerning density isoperimetric problems of the type
\[ \min \{ \Per_\nu(\Omega) :  \int_\Omega \nu(x) dx = c\}.\]
Note that here, the constant volume constraint is replaced by a constant density integral condition. This does not affect the proof of the $\Gamma$-convergence results, as long as $\nu$ is bounded on $D$. Indeed, the density constraints pass to the limit if $\nu \in L^\infty(D)$ and the construction of the recovery sequences in the proof of the (LS) property can be done in a similar way.

\begin{figure}
\centering
\includegraphics[width = 0.2 \textwidth]{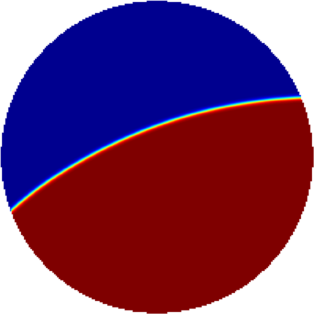}~
\includegraphics[width = 0.2 \textwidth]{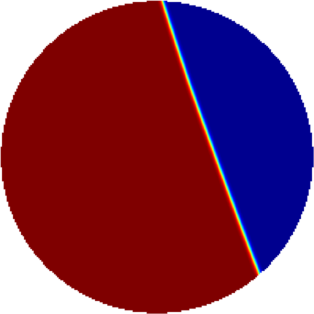}~
\includegraphics[width = 0.2 \textwidth]{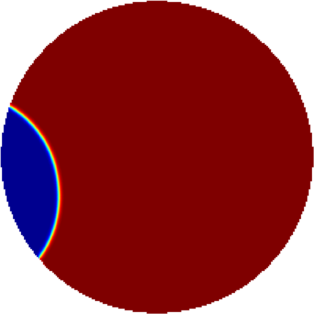}

\includegraphics[width = 0.2 \textwidth]{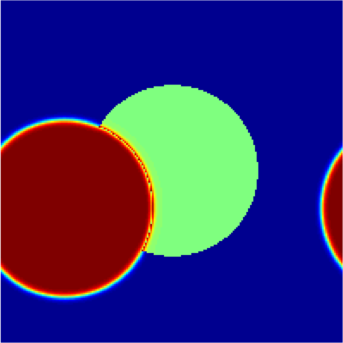}~
\includegraphics[width = 0.2 \textwidth]{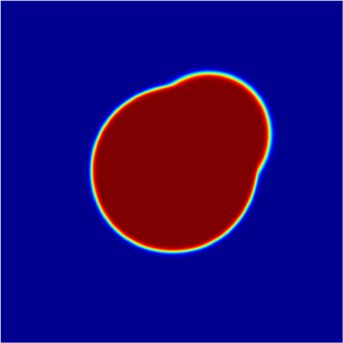}~
\includegraphics[width = 0.2 \textwidth]{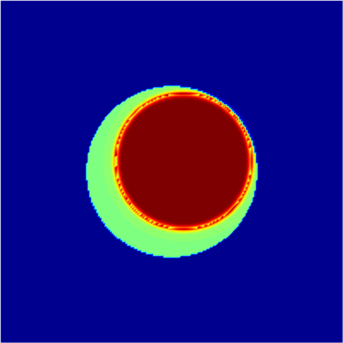}

\caption{Density equal to $\lambda$ on the unit disk and $1$ outside. First row $\lambda>1$, second row $\lambda<1$. }
\label{weighted1}
\end{figure}

We test the method in some of the cases presented in \cite{camivi} and \cite{morgan-pratelli}. The first case consists of a density which is equal to $\lambda$ inside the unit disk and $1$ outside. There are multiple situations concerning the parameter $\lambda$ and the fraction of the volume considered. We are able to recover numerically all the optimal shapes predicted theoretically. The results can be seen in Figure \ref{weighted1}

\bibliography{../master}
\bibliographystyle{plain}

\end{document}